\documentclass{amsart}
\usepackage{amscd, amssymb, portland, calc, ifthen}
\usepackage[matrix,arrow,curve,cmtip]{xy}
\CompileMatrices

\newtheoremstyle{mythm}{}{}%
  {\itshape}
  {}
  {\bfseries}
  {}
  { }
  {\thmnumber{#2.\hspace{1.5mm}}\thmname{#1}\thmnote{ #3}.}

\newtheoremstyle{myrmk}{}{}%
  {}
  {}
  {\bfseries}
  {}
  { }
  {\thmnumber{#2.\hspace{1.5mm}}\thmname{#1}\thmnote{ #3}.}

\numberwithin{equation}{subsection}

\theoremstyle{mythm}
\newtheorem{thm}[subsection]{Theorem}
\newtheorem{prop}[subsection]{Proposition}
\newtheorem{lemma}[subsection]{Lemma}

\theoremstyle{myrmk}
\newtheorem{rmk}[subsection]{Remark}

\newtheorem{questions}[subsection]{Questions}
\newtheorem{const}[subsection]{Construction}

\theoremstyle{plain}
\newtheorem*{thm*}{Theorem}
\newtheorem*{prop*}{Proposition}
\newtheorem*{cor*}{Corollary}

\theoremstyle{definition}
\newtheorem*{ack}{Acknowledgments}

\newcommand{\mathbold}{\bf}
\newcommand{\pser}[1]{{[\![{#1}]\!]}}

\newcommand{\pf}{^{\mathrm{pf}}}
\newcommand{\sep}{^{\mathrm{sep}}}

\newcommand{\red}[1]{{#1}_{\mathrm{red}}}

\newcommand{\longmap}{{\,\longrightarrow\,}}

\newcommand{\Hom}{{\mathrm{Hom}}}

\newcommand{\aut}{{\mathrm{Aut}}}

\newcommand{\Gal}{{\mathrm{Gal}}}

\newcommand{\gr}{{\mathrm{gr}}}

\newcommand{\ff}{{\mathbold F}}
\newcommand{\qq}{{\mathbold Q}}

\newcommand{\zz}{{\mathbold Z}}

\newcommand{\nn}{{\mathbold N}}

\newcommand{\texthidot}{{^{_{\bullet\!}}}}
\newcommand{\textlodot}{{_{^{\bullet\!}}}}

\newcommand{\m}{{\mathfrak{m}}}
\newcommand{\card}[1]{\,|#1|\,}
\newcommand{\codim}{{\mathrm{codim}\,}}
\newcommand{\Fr}{{\mathrm{Fr}}}

\newcommand{\comment}[1]{}

\newcommand{\un}{^{\mathrm{u}}}

\newcommand{\gc}{^{\mathrm{g}}}
\newcommand{\crp}[1]{{\mathsf{CRP}_{\negmedspace{#1}}}}

\newcommand{\arnv}[1]{{\mathrm{{ar_{\scriptstyle\mathrm{n}}^{#1}}}}}
\newcommand{\arn}{{\mathrm{{ar_{\scriptstyle\mathrm{n}}}}}}

\newcommand{\swk}{{\mathrm{sw_{\scriptscriptstyle K}}}}
\newcommand{\art}{{\mathrm{ar}}}

\newcommand{\p}{\mathfrak{p}}
\newcommand{\disc}{\mathfrak{D}}
\newcommand{\palg}[1]{\mathsf{PfAlg}_{{#1}}}

\newcommand{\barb}[1]{{\overline{{#1}}}}
\renewcommand{\m}{\p}
\newcommand{\qnn}{\qq_{\geq 0}}

\makeatletter
\def\@seccntformat#1{\@ifundefined{#1@cntformat}%
{\csname the#1\endcsname\quad}
{\csname #1@cntformat\endcsname}
}
\def\section@cntformat{\thesection.\enspace}
\def\subsection@cntformat{\thesubsection.}
\makeatother

\begin{document}

\newcommand{\kk}{{k}}
\newcommand{\kkt}{{\tilde{k}}}
\newcommand{\kkgc}{{k\gc}}

\newcommand\mnote[1]{}
\newcommand\bel[1]{{\mnote{#1}}\begin{equation}\label{#1}}
\newcommand\lb[1]{\label{#1}\mnote{#1}}

\newcounter{hour}\newcounter{minute}
\newcommand{\printtime}{\setcounter{hour}{\time/60}%
	\setcounter{minute}{\time-\value{hour}*60}%
	\ifthenelse{\value{hour}<10}{0}{}\thehour:%
	\ifthenelse{\value{minute}<10}{0}{}\theminute}

\setlength{\marginparwidth}{3cm}

\title{Conductors and the moduli of residual perfection}
\author[J.~Borger]{James M. Borger}
\address{University of Chicago, Chicago, Illinois, USA}
\email{borger@math.uchicago.edu}
\begin{abstract}
Let~$A$ be a complete discrete valuation ring with possibly imperfect
residue field.  The purpose of this paper is to give a notion of
conductor for Galois representations over~$A$ that generalizes the
classical Artin conductor.  The definition rests on two results of
perhaps wider interest: there is a moduli space that
parametrizes the ways of modifying~$A$ so that its residue field is
perfect, and any Galois-theoretic object over~$A$ can be recovered
from its pullback to the (residually perfect) discrete valuation ring
corresponding to the generic point of this moduli space.
\end{abstract}

\maketitle
\section*{Introduction}
Let~$A$ be a complete discrete valuation ring of residue
characteristic~$p>0$.  If the residue field of~$A$ is perfect, there
is a satisfactory theory~\cite[IV, VI]{Serre:CL} of ramification
over~$A$.  For example, let~$\rho$ be a
Galois representation over~$A$, which is
a continuous action of an absolute Galois group of the fraction field
of~$A$ on a finite-rank complex (for now) vector space.  Then there is a
non-negative integer, the Artin conductor of~$\rho$, that
measures the extent to which~$\rho$ is ramified. 
If, on the other hand, we allow the residue field of~$A$ to be
imperfect, ramification over~$A$ is still quite mysterious.
This prevents us from understanding, say,
ramification in codimension one of local systems on arithmetic surfaces.

The work in this paper began with the observation that much about
ramification over~$A$ can be understood by simply changing base to
sufficiently generic extensions with perfect residue field and of
relative ramification index one.  The first main point
(\ref{thm-moduli-rep}) is that such extensions make up the
(perfect-field-valued) points of a natural representable moduli
problem.  The {\em universal residual perfection}~$A\un$ of~$A$ is the
$A$-algebra corresponding to the representing object itself.  It is
not a discrete valuation ring---but, in a sense, only because its
residue ring is not a field (and actually not even noetherian).  The
{\em generic residual perfection}~$A\gc$ of~$A$ is the $A$-algebra
corresponding to the fraction field of the representing object.  It is
a complete discrete valuation ring with perfect residue field.
Both~$A\un$ and~$A\gc$ are, of course, unique up to unique
isomorphism.

It may be useful to keep a geometric analogue in mind. If we think of
complete residually perfect discrete valuation rings as being like
germs of curves, then it is natural to regard~$A\un$ as the universal
jet on~$A$ transverse to the maximal ideal and~$A\gc$ as the generic
jet.  (It would be interesting to see if the more sophisticated tools
from the theory of jet spaces have any relevance here.)

There is an explicit description
(\ref{thm-expl-Cfrp}) of these rings in the first section.  For example,
if~$A=\ff_p(x)\pser{y}$, then we have
\begin{align*}
A\un & \cong \ff_p(\bar{x})[u_1, u_2,\dots]^{p^{-\infty}}\pser{y} \text{\ and } \\
A\gc & \cong \ff_p(\bar{x},u_1, u_2,\dots)^{p^{-\infty}}\pser{y},
\end{align*}
where the~$A$-algebra structures are determined by the data~$x\mapsto
\bar{x}+u_1y+u_2y^2+\cdots$ and~$y\mapsto y$.
If~$A=\widehat{\zz[x]_{(p)}}$, we have
\begin{align*}
A\un & \cong W\bigl(\ff_p(\bar{x})[u_1, u_2,\dots]^{p^{-\infty}}\bigr) \text{\ and } \\
A\gc & \cong W\bigl(\ff_p(\bar{x},u_1, u_2,\dots)^{p^{-\infty}}\bigr),
\end{align*}
where~$W$ denotes the functor of Witt vectors
and the~$A$-algebra structures are determined
by~$x\mapsto(\bar{x},u_1^p,u_2^{p^2},\dots)$.

In the second section, I give some properties of~$A\un$
and~$A\gc$.  The most important is that the fraction field of~$A$
is algebraically closed in the fraction field of~$A\gc$.  And hence
the second main point: a Galois representation~$\rho$ over~$A$ is
determined by its pullback~$\rho|_{A\gc}$ to~$A\gc$.  Therefore, any
invariant of~$\rho$, such as a measure of ramification,
can be recovered (in principal) from~$\rho|_{A\gc}$.

Finally, I argue that in defining ``non-logarithmic'' conductors for
Galois representations over~$A$, the most simple-minded way of
proceeding along these lines is correct: the conductor~$\art(\rho)$
of~$\rho$ should be the classical Artin conductor of~$\rho|_{A\gc}$.
(Logarithmic conductors are more subtle.  See below.)
Taking this as the definition, we have the following result.

\begin{thm*}
Let~$\rho$ and~$\rho'$ be Galois representations over~$A$.
\begin{enumerate}
\item $\art(\rho)$ is a non-negative integer.
\item $\art(\rho\oplus\rho')=\art(\rho)+\art(\rho')$
\item If~$\rho$ is trivialized by a residually separable extension of~$A$,
  then~$\art(\rho)$ agrees with the classical Artin conductor of~$\rho$.
\item $\art(\rho)$ is zero if and only if~$\rho$ is unramified.
\item $\art(\rho)$ equals the codimension of the subspace of inertia
  invariants if and only if~$\rho$ is tame.
\end{enumerate}
\end{thm*}

This theorem is an elementary consequence of the basic properties
of~$A\gc$ proved in section~\ref{sec-prop}.  The proof of a slightly
stronger version is written down in section~\ref{sec-cond}.

For Galois representations of rank one, Kato~\cite{Kato:Imp} has
introduced a logarithmic conductor.  As first observed in the
work~\cite{Matsuda:Swan} of Matsuda (who credits the observation to
T. Saito), there is a non-logarithmic variant of Kato's conductor.  In
a companion paper, I show this non-logarithmic
conductor is invariant under pullback to~$A\gc$ and obtain the
following consequence.

\begin{cor*}[\cite{Borger:Kato-GRP}]
Let~$\rho$ be a rank-one Galois representation over~$A$.  Then the
non-logarithmic Kato conductor of~$\rho$ agrees with~$\art(\rho)$.
\end{cor*}

Kato's original, logarithmic conductor is not, however, always
invariant under pullback to~$A\gc$, and so the naive logarithmic
analogue of~$\art(\cdot)$ does not necessarily agree with Kato's
conductor.  For some brief thoughts on logarithmic
conductors for representations of higher rank, see~\ref{rmk-log-cond}.

Abbes and Saito~\cite{Abbes-Saito:Ram-I} have defined two $\qnn$-indexed
``upper'' filtrations of absolute Galois groups of the fraction fields
of complete discrete valuation rings with arbitrary residue fields.
At the end of section~\ref{sec-cond}, I give a definition of such a
filtration, which is compatible in the sense of breaks with this
paper's conductor.  It is tempting to hope it agrees with Abbes
and Saito's non-logarithmic filtration (shifted by one).

Boltje-Cram-Snaith~\cite{Boltje.et.al:Cond} and
Zhukov~\cite{Zhukov:Ram} also have approaches to non-abelian
ramification theory; the relations with them are even more mysterious.

\comment{
\marginpar{Add here.}
\begin{ack}
I would like to thank the Massachusetts Institute of Technology for
its hospitality and the National Science Foundation for support
while I was writing this paper.  I was motivated to think
about the material here by my thesis (Berkeley, 2000), which I
wrote under the advisement of Hendrik Lenstra, who also made several
useful comments on this paper.
\end{ack}
}

\setcounter{section}{0}
\section*{Conventions}
\lb{sec-notation}

All rings are commutative and contain 1, and all ring maps preserve 1.
The fraction field of a domain~$R$ is denoted~$\Fr(R)$.

An extension of a field is a homomorphism to another field.  An extension
of a discrete valuation ring is an injective local homomorphism to
another discrete valuation ring.  In both cases, we usually refer
to the target of the morphism, rather than the morphism itself, as the
extension.

Throughout,~$p$ denotes a fixed prime number, and $A$ is
a discrete valuation ring, held fixed in each subsection, whose
residue field has characteristic~$p$.  Its fraction field is
denoted by~$K$.  We say such a ring~$A$ is of {\em mixed
characteristic} if~$K$ has characteristic 0, and is of
{\em equal characteristic} if~$K$ has characteristic~$p$.  Variants of the
words {\em residue} and {\em generic} refer to the residue and
fraction fields of~$A$.  We also use the
same words to refer to extensions.
For example, we might say~$B/A$ is residually purely inseparable or is
generically Galois.  An extension~$B/A$ is {unramified}
(resp.\ {tame}) if it is residually separable and its ramification
index is one (resp.\ not divisible by~$p$).  
The notations~$e_{B/A}, f_{B/A},$ and~$f\sep_{B/A}$ denote
the ramification index, residue degree, and separable residue degree
of the extension~$B/A$.

\newcommand{\newpi}{{\tilde{\pi}}}
\section{The moduli of residual perfection}
\lb{sec-def}

The purpose of this section is to define the category~$\crp{A}$ of
complete residual perfections of~$A$, to prove the objects are
parametrized by a moduli space, and to give a concrete description of
this moduli space.  It is common in such matters that it is easier to
come up with the proofs than the statements.
It is no different in this section, and for this reason,
I have left a number of arguments to the reader.  Our general
reference for categorical terminology will be Mac Lane's
book~\cite{MacLane:CWM}.

\subsection{}
An~$\ff_p$-algebra~$R$ is {\em perfect} if the
endomorphism~$F\!:\!x\mapsto x^p$ of~$R$ is an isomorphism.  The
perfection~$R\pf=R^{p^{-\infty}}$ of an~$\ff_p$-algebra~$R$ is the universal
perfect~$\ff_p$-algebra that~$R$ maps to.  It is the colimit of the
iterates of~$F$.  For any~$\ff_p$-algebra~$S$, let~$\palg{S}$ be the
full subcategory of the
category of~$S$-algebras whose objects are perfect.

\subsection{}
Let~$\crp{A}$ be the full subcategory of the category of~$A$-algebras
consisting of objects~$B$ such that~$B$ is flat (i.e.,~$B$ is
torsion-free as an~$A$-module),~$B$ is complete with respect to the
ideal~$\p_A B$, and~$B/\p_AB$ is perfect.  (Note that the second
condition forces morphisms to be continuous.)  For an
object~$B\in\crp{A}$, let~$\p_B$ denote the ideal~$\p_AB$ (which might
not be prime), let~$\barb{B}$ denote $B/\p_B$, and if~$x$ is an
element of~$B$, let~$\bar{x}$ denote its image in~$\barb{B}$.
If~$f:B\to B'$ is a morphism in~$\crp{A}$, let~$\bar{f}$ denote its
reduction~$\barb{B}\to\barb{B'}$.
Let~$s_B$ denote the unique multiplicative section~\cite[II \S4 Prop.\
8]{Serre:CL} of the map~$B\to \barb{B}$.  Since~$s_B(x)$ is the unique
lift of~$x$ that has a~$p^m$-th root for all integers~$m$, every
morphism~$f:B\to B'$ in~$\crp{A}$ satisfies
\bel{eq-booha}
f\circ s_B=s_{B'}\circ\bar{f}.
\end{equation}

\subsection{}
\lb{sbsc-coeffs}
If~$\pi\in A$ is a uniformizer and~$x$ is an element of an
object~$B\in\crp{A}$, then there are unique
elements~$x_0,x_1,\dots\in\barb{B}$ such
that~$x=s_B(x_0)+s_B(x_1)\pi+\cdots$.  We call these the {\em
coefficients of~$x$ (with respect to~$\pi$)}.  If~$f:B\longmap C$ is a
morphism in~$\crp{A}$, then~(\ref{eq-booha}) implies that the
coefficients of~$f(x)$ are simply the images under~$\bar{f}$ of the
coefficients of~$x$.

\begin{thm}
\lb{thm-moduli-rep} 
The category~$\crp{A}$ has an initial object~$A\un$, and the functor
\[
\crp{A} \to \palg{\barb{A\un}}, \quad B \mapsto \barb{B}
\]
is an equivalence of categories.
\end{thm}

In other words, the functor~$\crp{A}\to\palg{\barb{A}}\:$ defined
by~$B\mapsto \barb{B}$ is a representable moduli problem with
universal object~$A\un$.


As usual, we can prove such a result in an abstract way or in a
concrete way.  I will give (or sketch) both.  The abstract proof
consists in showing directly that the
functor~$\crp{A}\to\palg{\barb{A}}$ is a cofibration in groupoids
(\ref{prop-cofib}) and then using Freyd's method to construct~$A\un$.
While appealing, the abstract proof has the big disadvantage that it
appears completely unreasonable to hope to use it to prove some
things---for example, that~$A\un$ is a domain.  The concrete proof, on
the other hand, consists in simply giving a presentation of~$A\un$ in
terms of a $p$-basis of~$\barb{A}$ lifted to~$A$.  The problem with
this approach is that it depends on choice, and the manner in which
the presentation of~$A\un$ depends on this choice is not clear.  It
would be nice to have an explicit description of~$A\un$ in terms of
familiar canonical constructions.

In any event, nothing that follows the abstract proof actually makes
use of it, and the reader so inclined can skip
to~\ref{prop-basic-structure} and take the canonical isomorphisms
in~\ref{sbsc-expl2} as definitions.

\begin{prop}
\lb{prop-cofib}
Let~$B$ be an object of~$\crp{A}$, and let~$\tilde{C}$ be a
perfect~$\barb{B}$-algebra.  Then there is an object~$C$ of~$\crp{A}$,
equipped with a morphism~$B\longmap C$ and a
morphism~$\tilde{C}\longmap\barb{C}$ of~$\barb{B}$-algebras, with the
property that for every other such object~$D$, there is a unique
map~$f:C\longmap D$ such that the diagram
\[\xymatrix{
{\barb{C}}\ar^{\bar{f}}[rr] && {\barb{D}} \\
  & {\tilde{C.}}\ar[ur]\ar[ul] \\}
\]
commutes. Moreover,~$C$ is unique up to unique isomorphism, and
the map~$\tilde{C}\to\barb{C}$ is an isomorphism.
\end{prop}

Intuitively,~$\tilde{C}$ has a lift to~$\crp{A}$ that is unique up to
unique isomorphism.

\begin{proof}
As always, the universality of~$C$ will determine it up to unique
isomorphism.  Let~$B[\tilde{C}]$ be the monoid algebra on the
multiplicative monoid underlying~$\tilde{C}$, and let~$C'$ be its
completion at~$\p_AB[\tilde{C}]$.  Because~$\tilde{C}$ is
perfect,~$C'$ is in~$\crp{A}$.  Now, let~$I$ be the kernel of the
obvious surjection~$\barb{C'}\longmap\tilde{C}$, let~$J$ be the ideal
of~$C'$ generated by the set~$s_{C'}(I)$, and put~$C=C'/J$.  The
map~$\tilde{C}\longmap\barb{C}$ is clearly an isomorphism.  The proof
that~$C$ satisfies the required universal property is left to the reader.
\end{proof}

\begin{lemma}
\lb{lemma-initial}
The category~$\crp{A}$ has an initial object.
\end{lemma}
\begin{proof}
We apply Freyd's method~\cite[V.6 Thm.\ 1]{MacLane:CWM}.  First, we
need to show~$\crp{A}$ has all (set-indexed) limits.  To do this, it
is enough to show it has arbitrary (set-indexed) products and
equalizers of pairs of arrows~\cite[V.2 Cor.\ 2]{MacLane:CWM}.

Fix a uniformizer~$\pi$ of~$A$.

Let~$f,g:B\longmap C$ be two morphisms in~$\crp{A}$.  Let~$E\subseteq B$
be their equalizer in the category of rings, and let~$\tilde{E}$
be the equalizer of~$\bar{f}$ and~$\bar{g}$ in the category of rings.
Then by~\ref{sbsc-coeffs}, we see that~$E$ consists of the elements
whose coefficients all lie in~$\tilde{E}$.
Therefore,~$\p_B^i\cap E = \p_A^i E$ for all~$i$.  It follows
that~$E$ is complete and the natural map~$E/p_AE\to\tilde{E}$ is an
isomorphism.  It is not hard to check~$\tilde{E}$ is perfect, and~$E$
is clearly torsion-free.  Hence~$E$ is an object of~$\crp{A}$.

The reader can check that~$\crp{A}$ has products.

We will now see that there exists a cardinal number~$\kappa$ such that
any object~$B$ in~$\crp{A}$ has a subobject whose cardinality is at
most~$\kappa$.  Let~$\tilde{C}$ be
the perfect subring of~$\barb{B}$ generated by the coefficients of the
images in~$B$ of elements of~$A$.  Then the cardinality of~$\tilde{C}$
is bounded (as~$B$ varies).  Let~$C$ be the closure of
the~$A$-subalgebra of~$B$ generated by the set~$s_B(\tilde{C})$.
Then the cardinalities of such subrings~$C$ are also bounded by some cardinal
number~$\kappa$.  It is easy to check that~$C/\pi C=\tilde{C}$,
and~$C$ is clearly complete and flat over~$A$.  Therefore,~$C$ is an
object of~$\crp{A}$.

Let~$D$ be the product of all objects~$B\in\crp{A}$ of cardinality at
most~$\kappa$.  Then by (the proof of) Freyd's theorem~\cite[V.6 Thm.\
1]{MacLane:CWM}, the equalizer of the set of endomorphisms of~$D$ is
an initial object.
\end{proof}

\begin{proof} (of~\ref{thm-moduli-rep})
Follows by general reasoning
from~\ref{prop-cofib} and~\ref{lemma-initial}.
\end{proof}

\subsection{}
The initial object $A\un$ of~$\crp{A}$ is called the {\em universal residual
perfection} of~$A$.  It is unique up to unique isomorphism.  Its
reduction~$\barb{A\un}$ (or, more properly, the spectrum of its
reduction) deserves to be called the moduli space of complete residual
perfections of~$A$, especially because its points with values in
perfect fields correspond to extensions of~$A$ of ramification index
one:

\begin{prop}
\lb{prop-basic-structure}
Let~$B$ be an object of~$\crp{A}$.  If~$\barb{B}$ is noetherian or a
domain, then~$B$ is the same.  If~$\barb{B}$ is a field, then~$B$ is a
discrete valuation ring.
\end{prop}

\begin{proof}
Suppose~$\barb{B}$ is noetherian.  Since~$\gr(B)$ is isomorphic
to~$\barb{B}[X]$, it is noetherian.  This
implies~$B$ is noetherian~\cite[10.25]{Atiyah-Macdonald:CA}.

Now suppose~$\barb{B}$ is a domain
and~$a$ and~$b$ are non-zero elements of~$B$ with~$ab=0$.
Since a uniformizer~$\pi$ of~$A$
is not a zero-divisor in~$B$, we can assume~$a,b\not\in\pi B$,
but this immediately
contradicts the fact that~$\barb{B}$ is a domain.

When~$\barb{B}$ is a field, it is easy to see that every element
outside~$\p_B$ is a unit.  Thus~$B$ is a noetherian local ring whose
maximal ideal is generated by a non-nilpotent element and, hence, a
discrete valuation ring~\cite[I \S 2]{Serre:CL}.
\end{proof}

\subsection{}
\lb{sbsc-expl}
Let us now introduce the notation and basic results we will need for the
explicit version of~\ref{thm-moduli-rep}.  

Let~$\pi\in A$ be a uniformizer, and let~$T$ be a lift to~$A$ of a
$p$-basis of~$\barb{A}$.  (A good general reference for information on
$p$-bases is
EGA~\cite[Ch.\ 0 \S 21]{EGA-no.20}.)  Let~$R_T$ be the polynomial
algebra~$\barb{A}[u_{t,j}\mid t\in T, j\in\zz_{>0}]$.  We will see
below that~$R_T\pf$ is naturally~$\barb{A\un}$, the moduli space we
seek.

\begin{lemma}
\lb{lemma-p-base-map}
Let~$Q$ be a residually perfect discrete valuation ring that is a
subring of~$A$ with the property that~$A/Q$ is an extension of
ramification index one.  Let~$B$ be a~$Q$-algebra that is complete
with respect to an ideal~$I$ that contains the image of the maximal
ideal of~$Q$.  Let~$n$ be a positive integer, and let~$\varphi':A\to
B/I^n$ be a~$Q$-linear homomorphism.  For every~$t\in T$, let~$x_t\in
B$ be a lift of~$\varphi'(t)$.  Then there is a unique $Q$-linear
map~$\varphi:A\to B$ such that~$\varphi'=\varphi\mod I^n$
and~$\varphi(t)=x_t$ for all~$t\in T$.
\end{lemma}
\begin{proof}
Since~$B$ is complete with respect to~$I$, it suffices by induction to
prove the existence and uniqueness of~$\varphi$ modulo~$I^{n+1}$.
Let~$\Omega^1_{A/Q}$ denote the~$A$-module of K\"ahler differentials
with respect to~$Q$.
Since~$A/Q$ is formally smooth~\cite[19.7.1]{EGA-no.20}, some
lift~$A\to B/I^{n+1}$ of~$\varphi'$ exists, and so the set of such
lifts is a torsor under
\begin{eqnarray*}
\Hom_A(\Omega^1_{A/Q},I^n/I^{n+1}) & = &
  \Hom_{\barb{A}}(\barb{A}\otimes_A \Omega^1_{A/Q},I^n/I^{n+1}) \\
 & = & \Hom_{\barb{A}}(\bigoplus_{t\in T} \barb{A}dt,I^n/I^{n+1})
\end{eqnarray*}
with the obvious action.
It follows that the image of~$T$ can be lifted arbitrarily and that
any such lift determines~$\varphi\!\!\mod I^{n+1}$.
\end{proof}

\begin{const}
\lb{const-UV}
{\em Functors~$U:\crp{A}\to\palg{R_T\pf}$
and~$V:\palg{R_T\pf}\to\crp{A}$.}

Let~$B$ be an object of~$\crp{A}$.  For~$t\in T$ and~$j\in\zz_{>0}$,
let~$v_{t,j}\in\barb{B}$ denote the~$j$-th coefficient (\ref{sbsc-coeffs})
with respect to~$\pi$ of the image of~$t$ in~$B$.
The data~$u_{t,j}\mapsto v_{t,j}$
gives~$\barb{B}$ the structure of an~$R_T$-algebra.  Since~$\barb{B}$ is
perfect, it has a unique compatible $R_T\pf$-algebra
structure;~$U(B)$ is then~$\barb{B}$ with this~$R_T\pf$-algebra structure.
It is easy to see this is functorial. 



Let~$S$ be a perfect~$R_T\pf$-algebra.
If~$A$ is of equal characteristic, set~$V(S) = S[\![{\newpi}]\!]$,
where~$\newpi$ is a free variable.  Then~$V(S)$ is
an~$\ff_p\pser{{\newpi}}$-algebra, and~${\newpi}\mapsto
\pi$ makes~$A$ into an~$\ff_p\pser{{\newpi}}$-algebra.
By~\ref{lemma-p-base-map}, there is a
unique~$\ff_p\pser{\newpi}$-linear map~$A\to V(S)$ such that for
all~$t\in T$,
\[
t\mapsto\bar{t}+\sum_{j\geq 1}u_{t,j}\newpi^j.
\]

If~$A$ is of mixed characteristic, let~$C$ be the Cohen
subring~\cite[pp.\ 82-83]{Cohen:CLR} of~$A$ determined by~$T$.
It is a complete discrete valuation ring of absolute
ramification index one that contains~$T$, and~$A/C$ is a finite
residually trivial extension.
Let~$W$ be the ring of Witt vectors~\cite[0.1]{Illusie:dRW} with
coefficients in~$S$ and let~$s_W:S\to W$ denote the Teichm\"uller
section.  Again by~\ref{lemma-p-base-map} (taking~$Q=\zz_p$),
there is a unique map~$C[X]\to
W[\![\newpi]\!]$ such that~$X\mapsto \newpi$ and for all~$t\in T$, we have
\[
t\mapsto s_{W}(\bar{t})+\sum_{j\geq 1}s_{W}(u_{t,j})\newpi^j.
\]
View~$A$ as a quotient of~$C[X]$ using the map~$X\mapsto \pi$, and
put~$V(S) = A \otimes_{C[X]} W\pser{\newpi}$.

In either case, it is easy to see that~$V$ is a functor 
from~$\palg{R\pf_T}$ to the category of~$A$-algebras.
\end{const}

\begin{prop}
The image of~$V$ is in~$\crp{A}$.
\end{prop}
\begin{proof}
Let~$S$ be a perfect~$R_T\pf$-algebra.  Of the properties~$V(S)$ is
required to satisfy to be in~$\crp{A}$, the only one that is not
immediately clear is flatness over~$A$ in mixed characteristic.
To show this, it suffices to show the element~$\pi\otimes 1 = 1\otimes \newpi$
is not a zero-divisor in~$V(S)$.

Let~$g(X)$ denote the Eisenstein polynomial that generates the kernel
of the surjection~$C[X]\to A$,
and let~$h(\newpi)$ denote its image in~$W\pser{\newpi}$.  We will
show~$\newpi$ is not a zero-divisor in the
ring~$V(S)=W\pser{\newpi}/(h(\newpi))$.
Let~$f_1(\newpi)$ and $f_2(\newpi)$ be elements of~$W\pser{\newpi}$ such
that~$\newpi f_2(\newpi)=h(\newpi)f_1(\newpi)$.

Suppose~$f_1(0)\neq 0$.  Then~$h(0)$, which is the image of~$g(0)$
under the map~$C\to W$, is a zero-divisor.  But since~$g(X)$ is an
Eisenstein polynomial, this implies~$p$ is a zero-divisor in~$W$,
which is impossible.  Therefore, we have~$f_1(0)=0$ and, hence,
\[
f_2(\newpi) = h(\newpi)\frac{f_1(\newpi)}{\newpi} \in h(\newpi))W\pser{\newpi}.
\]
Thus~$f_2(\newpi)$ reduces to zero in~$V(S)$, and so~$\newpi$ is not a
zero-divisor in~$V(S)$.
\end{proof}

\begin{const}
{\em Natural transformations~$\eta:1\to UV$
and~$\varepsilon:VU\to 1$}

Let~$B$ be an object of~$\crp{A}$ and~$S$ be an object
of~$\palg{R_T\pf}$.  Let~$\eta(S)$ be the map
\[
S\longmap UV(S),\ \ x\mapsto\barb{s_{V(S)}(x)}.
\]
If~$A$ is of equal characteristic, let~$\varepsilon(B)$ be the
composite~$VU(B)=\barb{B}\pser{\newpi} \longmap B$ defined
by~$\newpi\mapsto\pi$ and~$b\mapsto s_B(b)$ for~$b\in\barb{B}$.  It is
a homomorphism of rings because~$s_B$ is~\cite[II Prop.\ 8]{Serre:CL}.

If~$A$ is of mixed characteristic, there is a unique
map~$W(\barb{B})\longmap B$ that reduces to the
identity~\cite[II Prop.\ 10]{Serre:CL}.
Then,~$\newpi\mapsto \pi$ determines a
map~$W(\barb{B})\pser{\newpi}\longmap B$ and, hence, a map
\[
\varepsilon(B):VU(B) = V(\barb{B}) = A\otimes_{C[X]}W(\barb{B})\pser{\newpi} \longmap B.
\]
It is easy to see~$\varepsilon$ is natural in~$B$ and~$\eta$ is natural
in~$S$.
\end{const}

\begin{thm}
\lb{thm-expl-Cfrp}
$\langle V,U;\eta,\varepsilon\rangle$ is an adjoint equivalence
between~$\crp{A}$ and~$\palg{R_T\pf}$.
\end{thm}

The proof is nothing more than a straight-forward verification of the
so-called triangular identities ($\varepsilon V \circ V \eta = 1$
and~$U\varepsilon \circ \eta U = 1$) and is left to the reader.

\subsection{}
\lb{sbsc-expl2}
It follows that the unique morphism~$A\un\to V(R_T\pf)$ is an
isomorphism.  This, in turn, induces a canonical
isomorphism~$\barb{A\un}=R_T\pf$, and so~$A\un$ is,
by~\ref{prop-basic-structure}, a domain.  ($A\un$ is also
integrally closed, but we will not use this fact here.) The {\em
generic residual perfection}~$A\gc$ of~$A$ is the object corresponding
under~\ref{thm-moduli-rep} to the fraction field of~$\barb{A\un}$,
that is, the object of~$\crp{A}$ corresponding to the generic point of
the moduli space.  It is canonically isomorphic to~$V(\Fr(R_T\pf))$.
By~\ref{prop-basic-structure}, it is a complete discrete valuation
ring with perfect residue field.

\subsection{}
\lb{sbsc-Aun-func}
Note that if~$B$ is an extension of~$A$ of ramification index one,
then~$\crp{B}$ is a subcategory of~$\crp{A}$, and so there is a unique
map~$A\un\to B\un$. If it happens that~$T$ can be extended to a
lift~$T'\subset B$ of a $p$-basis for~$\barb{B}$, this map
is the same as the map associated by~\ref{thm-expl-Cfrp} to the
inclusion~$T\to T'$.  The functoriality of the generic residual
perfection is a little more subtle.  We will consider it in the next
section.

\subsection{} \lb{sbsc-filt}
It is not hard to show that the filtration
\begin{equation*}
F_nR_T\pf = \barb{A}[u_{t,j} \mid t\in T, 1\leq j\leq n]\pf
\end{equation*}
of~$R_T$ is independent of the choices of~$T$ and~$\pi$.  Though its
role in this paper is small, this filtration is important and
should not be ignored.  (See~\ref{rmk-log-cond}.)

\begin{rmk}
I have written this section in the context needed for this paper,
but I believe similar constructions could be made in much greater generality.

For example, it appears that the functor~$V$ could be defined for
rings~$A$ (or schemes, or \dots) much more general than discrete valuation
rings.  It is not yet clear, however, what category should
replace~$\crp{A}$ and what kind of rings~$A$ should be considered.
Also, even if we look only at discrete valuation rings~$A$, it is probably
more natural, from the abstract perspective,
to require that each object of~$\crp{A}$ be equipped with
a multiplicative section from its residue ring (and to require that
morphisms respect this structure) than to require that it be complete
and its residue ring be perfect.  

\end{rmk}

\section{Properties}
\lb{sec-prop}

As in the previous section,~$A\un$ and~$A\gc$ will denote the
universal and generic residual perfections of~$A$, and~$K\gc$ will
denote the fraction field of~$A\gc$.

\begin{prop}
\lb{prop-gc-func}
Let~$B$ be an extension of~$A$.  If~$e_{B/A}=1$
and~$B/A$ is residually separable,
then there exists a unique map~$A\gc\longmap B\gc$
of~$A$-algebras.
\end{prop}

Note that we do not require that~$B/A$ be residually algebraic, and
so~$B$ satisfies these assumptions
if and only if~$B/A$ is formally smooth~\cite[19.6.1,19.7.1]{EGA-no.20}.

\begin{proof}
By~\ref{thm-moduli-rep}, it is enough to show there is a unique
map~$\barb{A\gc}\longmap\barb{B\gc}$ of~$\barb{A\un}$ algebras.
Since~$\barb{A\gc}$ is the fraction field of~$\barb{A\un}$, there is
clearly at most one.  To show there is at least one, we just have to
check that the map~$\barb{A\un}\to\barb{B\gc}$ is an inclusion.  
Since any $p$-basis for~$\barb{A}$ can be extended to one for~$\barb{B}$,
this follows immediately from~\ref{thm-expl-Cfrp} and~\ref{sbsc-Aun-func}.
\end{proof}

\begin{prop}
\lb{prop-gen-etale-bc}
Let~$B$ be a finite \'etale extension of~$A$.  If~$C$ is an object
of~$\crp{A}$, then~$C\otimes_A B$ is an object of~$\crp{B}$.
\end{prop}
\begin{proof}
Since~$B\otimes_A\overline{C}$ is finite
\'etale over $\overline{C}$, which is perfect, it is perfect
(by, for example,~\cite[21.1.7]{EGA-no.20}).  Since~$C\otimes_A B$
is a finitely generated free~$C$-module, it is complete.
Therefore,~$C\otimes_A B$ is an object of~$\crp{B}$.
\end{proof}

\begin{prop}
\lb{prop-etale-bc}
Let~$B$ be a finite \'etale extension of~$A$.  Then the induced
maps
\begin{equation*}
B\otimes_A A\un\rightarrow B\un \quad\text{and}\quad
B\otimes_A A\gc\rightarrow B\gc 
\end{equation*}
of $B$-algebras are isomorphisms.
\end{prop}
\begin{proof}
By~\ref{prop-gen-etale-bc}, both~$B\otimes_A A\un$
and~$B\otimes_A A\gc$ are in~$\crp{B}$, and so
by~\ref{thm-moduli-rep}, it is enough to show they are isomorphisms
after tensoring with~$\barb{B}$.  Since~$\barb{B}/\barb{A}$ is finite
and separable, any $p$-basis for~$\barb{A}$ is one for~$\barb{B}$.
Again,~\ref{thm-expl-Cfrp} and~\ref{sbsc-Aun-func} complete the proof.
\end{proof}

The rest of the results in this section are devoted to the proof of the
following theorem.

\begin{thm}
\lb{thm-Galois-surj}
Fix a separable closure of~$K\gc$.  Then the map~$G_{K\gc}\to G_K$
of absolute Galois groups is surjective.  The induced maps of
inertia groups and wild inertia groups are also surjective.
\end{thm}

\begin{rmk}
It would be interesting to see if the analogous result is true for
some motivic Galois group.  For example, when~$A$ is of equal
characteristic, is the functor from crystals (of whatever kind) on~$K$
to crystals on~$K\gc$ fully faithful?
\end{rmk}

\begin{lemma}
\lb{lemma-alg-closed}
$K$ is algebraically closed in~$K\gc$.
\end{lemma}
\begin{proof}
If~$A$ is residually perfect, then~$K\gc=K$ and the result is trivially true.
Now assume~$A$ is not residually perfect.
Let~$L/K$ be a finite subextension of~$K\gc/K$, and
let~$B$ be the normalization of~$A$ in~$L$.  Since~$e_{A\gc/A}=1$, we
have~$e_{B/A}=1$.  Since~$\barb{A}$ is separably closed
in~$\overline{A\gc}$, the extension~$B/A$ is residually purely
inseparable.  It is therefore enough to show it is residually separable.

Suppose there is an element~$\bar{a}\in\barb{A}-(\barb{A})^p$ that has a
$p$-th root in~$\barb{B}$.  We will derive a contradiction (working
only modulo~$\m_A^2$).  Let~$a\in A/\m_A^2$ be a lift of~$\bar{a}$ and
let~$x\in B/\m_B^2$ be a lift of~$\sqrt[p]{\bar{a}}$.  Then there is
an element~$y\in\barb{B}$ such that~$a=x^p+\pi y$.

Let~$\theta$ denote the map~$B\longmap A\gc/\m_A^2 A\gc$, and let~$s$
denote the multiplicative section of the map~$A\gc/\m_A^2 A\gc
\longmap \overline{A\gc}$. Since the set~$\{\bar{a}\}$ can be extended to a
$p$-basis of~$\barb{A}$,
there is (by, say,~\ref{thm-expl-Cfrp}) an
element~$u\in\overline{A\gc}-\barb{A}\pf$ such that~$\theta(a)=
s(\bar{a}) + \pi u$.  Take~$v\in \barb{A\gc}$ such that~$\theta(x)=
s(\bar{x}) + \pi v$.  Then we have
\begin{eqnarray*}
s(\bar{a})+\pi u & = & \theta(a) \\
 & = & \theta(x)^p+\pi y \\
 & = & (s(\bar{x})+\pi v)^p + \pi y \\
 & = & s(\bar{x}^p)+\pi y 
\end{eqnarray*}
and thus~$u=y\in\barb{A}\pf$.  Since this cannot be, we have our contradiction.
\end{proof}

\subsection{}
An extension~$B$ of~$A$ is said to be {\em monogenic} if there is an
element~$x\in B$ that generates~$B$ as an~$A$-algebra.  For example,
any finite extension that is generically and residually separable is
monogenic~\cite[III \S 6 Prop.\ 12]{Serre:CL}.

\begin{lemma}
\lb{lemma-monogenic} 
Let~$B$ be a generically separable extension of~$A$.  Then~$B/A$ is
monogenic if and only if~$B\otimes_A A\gc$ is a discrete valuation
ring.  In this case, we have~$f\sep_{B\otimes_A A\gc/A\gc} =
f\sep_{B/A}$.
\end{lemma}
\begin{proof}
Suppose~$B\otimes_A A\gc$ is a discrete valuation ring.  Then
we have
\[
A\gc\otimes_A\Omega^2_{B/A}=\Omega^2_{B\otimes_A A\gc/A\gc}=0,
\]
where, as usual,~$\Omega^2_{*/*}$ denotes the second exterior power of
the module of relative K\"ahler differentials.  Therefore,~$\Omega^2_{B/A}$ is
zero.  The extension~$B/A$ is then monogenic by de
Smit~\cite[4.2]{DeSmit:Imp}.

Now suppose~$B/A$ is monogenic.
By~\ref{lemma-alg-closed}, the ring~$B\otimes_A A\gc$ is a
domain.  By~\ref{prop-etale-bc}, it suffices to
assume~$B/A$ is residually purely inseparable, and so it is enough
to show that~$B\otimes_A A\gc$ is generated as an~$A\gc$-algebra by the
root of an Eisenstein polynomial.

There is~\cite[Prop.\ 1]{Borger:Monogenic} an integral extension~$A'$
of~$A$ such that~$e_{A'/A}=f\sep_{A'/A}=1$ and~$B\otimes_A
A'$ is a discrete valuation ring with~$f_{B\otimes_A A'/A'}=1$.
By Zorn's lemma, we can assume the residue field of~$A'$
is~$\barb{A}\pf$.  Now let~$x\in B$ be a generator of
the~$A$-algebra~$B$ and let~$f(X)\in A[X]$ be its characteristic
polynomial over~$A$.  Put
\[
y = x\otimes 1 - 1\otimes s_{A\un}(\bar{x}) \in B\otimes_A A\un.
\]
Then~$y$ generates~$B\otimes_A A\un$ as an~$A\un$-algebra, the
element~$\bar{y}\in\barb{B\otimes_A A\un}=\barb{A\un}$ is zero, and
the polynomial
\[
g(X)=f(X+s_{A\un}(\bar{x}))\in A\un[X]
\]
is the characteristic polynomial of~$y$ over~$A\un$.  We will
show~$g(X)|_{A\gc}$ is an Eisenstein polynomial.

First, note that since~$\bar{y}=0$, we have~$g(X) \equiv X^n \mod \m_A A\un$,
where~$n$ is the degree of~$B/A$.
All that remains is to
show~$g(0)|_{A\gc} \notin \m_{A\gc}^2$.  Hence, it is enough to
show~$g(0)|_{A\un} \notin \m_A^2 A\un$ and, therefore,
even~$g(0)|_{A'} \notin \m_{A'}^2$.  But since~$y|_{A'}$
generates~$B\otimes_A A'$ as an~$A'$-algebra and since~$B\otimes_A A'$
is a discrete valuation ring,~$g(X)|_{A'}$ is an Eisenstein
polynomial.  Therefore,~$g(0)|_{A'} \notin \m_{A'}^2$ and~$g(X)|_{A\gc}$
is an Eisenstein polynomial.
\end{proof}

\begin{lemma}
\lb{lemma-inertia-surj}
Let~$B$ be a finite, generically Galois extension of~$A$.
Let~$B'$ be the integral closure of the domain~$B\otimes_A A\gc$.
Then~$f\sep_{B'/A\gc} = f\sep_{B/A}$.
\end{lemma}
\begin{proof}
By~\ref{prop-etale-bc}, it is enough to assume~$B/A$
is residually purely inseparable.
Let~$G$ be the generic Galois group of~$B/A$.
Then~\ref{lemma-alg-closed} implies~$G$ is also the generic Galois
group of~$B'/A\gc$.  Let~$C'$ be the maximal \'etale
subextension of~$B'/A\gc$.  Then~$C'$ corresponds to a normal
subgroup of~$G$ and, hence, to a generically
Galois subextension~$C$ of~$B/A$.  Since the extension~$B/A$ is
residually purely inseparable and $e_{C'/A}=1$, the ramification index
of~$C/A$ is one and its Galois group is a $p$-group.

Now let~$D$ be a monogenic subextension of~$C/A$.  Then~$D\otimes_A
A\gc$ is a discrete valuation ring by~\ref{lemma-monogenic}.
Since~$D\otimes_A A\gc/A\gc$ is a subextension of~$C'/A\gc$, it is \'etale.
Now, the residue field of~$D\otimes_A A\gc$ 
is~$\red{(\barb{D} \otimes_{\barb{A}} \overline{A\gc})}$.
Since~$\barb{D}$ is a finite purely inseparable extension of~$\barb{A}$,
this reduced quotient is~$\overline{A\gc}$.  Thus, the \`etale
extension~$D\otimes_A A\gc/A\gc$ is trivial, and therefore, so
is~$D/A$.

But since the generic Galois group of~$C/A$ is a $p$-group and since
all extensions of degree~$p$ are monogenic, the only way for all
monogenic subextensions of~$C/A$ to be trivial is if~$C/A$ itself is
trivial.  And this can happen only if~$C'/A\gc$ is trivial.
\end{proof}

\begin{proof}
(of~\ref{thm-Galois-surj})
It follows from~\ref{lemma-alg-closed} that~$K$ is
separably closed in~$K\gc$.  The surjectivity of~$G_{K\gc}\to G_K$
is just the translation of this into Galois theory.

The image of the inertia
subgroup of~$G_{K\gc}$ is then a closed normal subgroup~$N$ of~$G_K$
and is contained in the inertia subgroup of~$G_K$.
Let~$M$ be the corresponding extension of~$K$.  Now, for any finite
extension~$B$ of~$A$ that is contained in~$M$,
we see, by lemma~\ref{lemma-inertia-surj}, that~$B/A$ is unramified.
Therefore,~$N$ is the entire inertia group of~$K$.

Since the wild inertia groups are the unique (pro-)$p$-Sylow subgroups
of the corresponding inertia groups~\cite[Ch.\ I, 1.1]{Serre:Gal-coh},
the surjectivity of the map between wild inertia groups follows immediately.
\end{proof}

\begin{questions}
The maps~$A\to A\un$ and~$A\to A\gc$ are faithfully flat.
How can we best understand the descent data?  The generic
descent data?  What can be said about the structure of the exact sequence
\[
1 \to H \to G_{K\gc} \to G_K \to 1?
\]
Since~$A\gc$ is residually perfect, the classfield theory of
Hazewinkel and Serre~\cite{Hazewinkel:LCFT} gives a description of the
abelianization of~$G_{K\gc}$.  Does this give a
classfield theory describing the abelianization of~$G_K$?
When~$K$ is a higher local field, how does the classfield theory
of~$K\gc$ relate to Kato and Parshin's classfield
theory~\cite{Kato:LCFTK-I,Parshin:LCFT} of~$K$?
\end{questions}

\section{Conductors}
\lb{sec-cond}

\subsection{}
\lb{sbsc-cond-def}
Fix a field~$\Lambda$ whose characteristic is not~$p$, and
let~$\rho$ be a Galois representation over~$A$, that is, a
homomorphism~$\rho:G\longmap\aut(V)$, where~$G$ is the Galois group of a
finite generically Galois extension~$B$ of~$A$
and~$V$ is a finite-rank~$\Lambda$-module.  For any~$i\in\nn$,
let~$G_i$ be the kernel of the map~$G\to\aut(B/\m_B^{i+1})$.  Define
\[
\arnv{B}(\rho)={e^{-1}_{B/A}}\sum_{i\geq 0}\card{G_i}\codim V^{G_i},
\]
where~$\card{\cdot}$ denotes cardinality and~$\codim V^{G_i}$ is the
codimension of the subspace of~$G_i$-invariants.
We call~$\arnv{B}(\rho)$ the {\em naive Artin conductor of~$\rho$ with
respect to~$B$}.  It is a non-negative rational number.
If~$B/A$ is residually separable, then~$\arnv{B}(\rho)$ is left unchanged if
we replace~$B$ with a larger generically Galois and residually separable
extension of~$A$.  In this case, we will use the notation~$\arnv{}(\rho)$.

Let~$L$ denote the fraction field of~$B$. Then
by~\ref{thm-Galois-surj}, $L\otimes_K K\gc$ is a finite Galois
extension of~$K\gc$, and its Galois group is canonically isomorphic
to~$G$.  Let~$\rho|_{A\gc}$ denote the resulting representation
of~$\Gal(L\otimes_K K\gc/K)$. Define the Artin conductor~$\art(\rho)$
of~$\rho$ to be~$\arn(\rho|_{A\gc})$.  We have the following slightly
stronger version of the theorem in the introduction:

\begin{thm}
\lb{thm-properties}
Let~$A,\rho,G,V,$ and~$B$ be as above, and let~$\rho'$ be another
Galois representation over~$A$.  Then we have the following:
\begin{enumerate}
  \item $\art(\rho\oplus\rho') = \art(\rho)+\art(\rho')$.
         \lb{thm-properties-1}
  \item $\art(\rho)$ is a non-negative integer.
         \lb{thm-properties-2}
  \item $\art(\rho)=0$ if and only if~$\rho$ is unramified, i.e.,~$G_0$
         acts trivially on~$V$.  \lb{thm-properties-3}
  \item If~$B/A$ is monogenic, then~$\art(\rho) = \arnv{B}(\rho)$.
         \lb{thm-properties-5}
  \item The following are equivalent: \lb{thm-properties-4}
      \begin{enumerate}
        \item $\rho$ is tame, i.e., the $p$-Sylow subgroup of~$G_0$
              acts trivially on~$V$ \lb{thm-properties-4a} 
        \item $\art(\rho)=\codim(V^{G_0})$ \lb{thm-properties-4b}
        \item $\art(\rho)\leq \codim(V^{G_0})$ \lb{thm-properties-4c}
      \end{enumerate}
\end{enumerate}
\end{thm}
\begin{proof}
Statement~\ref{thm-properties-1} follows from the additivity of the
naive Artin conductor.  Statement~\ref{thm-properties-2} follows from the
classical Hasse-Arf theorem~\cite[VI \S 2]{Serre:CL} (if the
characteristic of~$\Lambda$ is not zero, see~\cite[19.3]{Serre:RLGF})
and the non-negativity of the naive Artin conductor.
Statement~\ref{thm-properties-3} follows from the statement about inertia
groups in~\ref{thm-Galois-surj}.

Now consider statement~\ref{thm-properties-5}.
Let~$G'$ be the generic Galois
group of the extension~$B\otimes_A A\gc/A\gc$. 
By~\ref{lemma-monogenic}, the tensor product~$B\otimes_A A\gc$
is a discrete valuation ring.   
A short argument then shows that for all~$i\in\nn$, we have
\[
G_i=G'_{ei}=\cdots =G'_{ei+e-1},
\]
where~$e=e_{B\otimes_A A\gc/B}$.  Therefore,~$\arnv{B}(\rho) =
\arnv{B\otimes_A A\gc}(\rho|_{A\gc}) = \art(\rho)$.

As for statement~\ref{thm-properties-4},
if~$\rho$ is tame, then there is some generically Galois and
residually separable extension~$B$ of~$A$ such that~$\rho|_B$ is
trivial.  Since~$B/A$ is residually separable, it is monogenic, and
so statement~\ref{thm-properties-5} implies~$\art(\rho)=m$.
If, on the other hand, we have~$\art(\rho)\leq m$, then~$\rho|_{A\gc}$ is
tame,~\ref{thm-Galois-surj} therefore implies~$\rho$ is tame.
\end{proof}

\begin{rmk}
\lb{rmk-log-cond}
As mentioned in the introduction, logarithmic conductors are more
subtle than non-logarithmic ones.  For~$\chi\in H^1(K,\qq/\zz)$,
let~$\swk(\chi)$ be Kato's Swan conductor~\cite[2.3]{Kato:Imp}, which
equals the logarithmic order of the pole of his refined Swan
conductor.  Then there can be classes~$\chi$ (those that
Matsuda says are {\em of type II}~\cite[3.2.10]{Matsuda:Swan}) such
that~$\swk(\chi)=\swk(\chi|_{A\gc})+1$.  Because of this, the naive
Swan analogue of~$\art(\cdot)$ does not always agree with Kato's Swan
conductor.  I believe that by taking into account the
filtration~$F_{\textlodot}\barb{A\un}$ of~$A\un$
(see~\ref{sbsc-filt}), one could give a good definition of a
logarithmic conductor.  Indeed, using the techniques in
Matsuda~\cite{Borger:Kato-GRP,Matsuda:Swan}, it is easy to see how to
recover the Kato-Swan conductor of~$\chi$ from the refined Swan
conductor of~$\chi|_{A\gc}$ in equal characteristic.
\end{rmk}
\begin{rmk}
It may also be worth mentioning that the conductor~$\art(\cdot)$ probably does 
not satisfy an induction formula.  Let~$A'/A$ be a finite generically
separable extension,~$\rho'$ a Galois representation over~$A'$, and~$\rho$
the induced representation over~$A$.
If~$\barb{A}$ is perfect, then~\cite[VI \S 2]{Serre:CL}
\bel{eq-ind}
\art(\rho) = f_{A'/A}\art(\rho) + \dim(\rho)v_A(\disc),
\end{equation}
where~$\disc$ denotes the discriminant of~$A'/A$.
In general, the refined conductor of~$\rho$ should be the norm, in
some suitable sense, of the refined conductor of~$\rho'$.
When~$\barb{A}$ is perfect, the conductor determines the refined
conductor up to an element of~$A^*$, and then~\ref{eq-ind}
would follow from such a norm formula for refined conductors.  But
when~$\barb{A}$ is not perfect, the refined conductor contains more
information than the conductor together with a unit.  In fact, even
in the abelian case, the
conductor of~$\rho'$ probably does even determine the conductor
of~$\rho$ in general.

\subsection{}
There is, however, an apparently satisfactory theory of the upper
filtration.  In the notation of~\ref{sbsc-cond-def}, since~$G$ is
naturally~$\Gal(L\otimes_K K\gc/K\gc)$, it inherits a~$\qnn$-indexed
upper filtration~\cite[IV \S3]{Serre:CL}.  This filtration is stable
under passage to quotients for the simple reason that the same is true
for residually separable extensions.  We also have the usual relation
to conductors: Following Katz~\cite[1.1]{Katz:Gauss}, there is a break
decomposition~$V=\oplus_x V(x)$, with respect to~$G^{\texthidot}$, and
we have
\[
\art(\rho) = \codim(V^{G_0}) + \sum_{x\in\qnn}x\dim V(x).
\]
\end{rmk}

\bibliography{references}
\bibliographystyle{plain}
\end{document}